# Production and distribution planning, scheduling, and routing optimization in a yogurt supply chain under demand uncertainty: A case study

Babak Javadi[1*], Zeinab Salimzadeh[1], Amir Hossein Akbari[2], Mahla Yadegari [1], Mohammadreza Abdali[1]

*[1] Department of Industrial Engineering, Faculty of Engineering, College of Farabi, University of Tehran, Iran
[2] Faculty of Industrial Engineering, Iran University of Science and Technology, Tehran, Iran

*Abstract:*

Considering the evolution of the food industry and its challenges, like high perishability, managing the food industry supply chain is a key focus for researchers and decision-makers. Uncertainty in decision-making has gained importance, particularly in the yogurt industry, known for its complexity. This study addresses production and distribution planning, scheduling, and routing in the yogurt supply chain. The problem is characterized by multiple products, a single plant, multiple distribution centers, multiple periods, and various transportation methods. A mixed-integer non-linear programming (MINLP) model is used to minimize total costs, including production, setup, overtime, unmet demand, and transportation. Additionally, a robust fuzzy programming approach is applied under uncertainty, with linearization procedures proposed to convert it into a linearized mixed-integer programming formulation. The problem is tested with two data types: a sample problem in three sizes (small, medium, and large) and real data from Kalle Dairy Company, Iran. A Genetic Algorithm (GA) is developed to solve the problem, with necessary modifications made for its application. The GA's performance is compared to an exact algorithm (Branch & Cut), showing that the company's production policy adapts daily to meet demand precisely. The shift to smaller batch production and longer shelf life allows better stock allocation and avoids shortages in uncertain conditions. The company's policies adapt to severe fluctuations in the business environment, though this requires high costs, such as inventory maintenance.

---

[*] Corresponding author: babakjavadi@ut.ac.ir(babakjavadi)

# 1. INTRODUCTION

The supply chain of perishable goods has always been one of the most critical and challenging issues in supply chain management at different times. However, many significant problems must be addressed in the food supply chain, such as high chain costs, food waste, and environmental issues. The ever-rising expansion of food production networks, as well as the diversification of producers and the various ways in which food products are supplied, and on the other hand, the existence of irreversible uncertainties, such as uncertainty in market demand and prices, complicates the management of the food supply chain.

Within the dairy industry, production planning of yogurt products is the most challenging task since we are facing vast variants of products. M. Lutke Entrup (1) considered shelf life the most distinctive factor in fresh food production planning. They mentioned that shelf-life restrictions directly influence wastage, out-of-stock rate, and inventory level. Amorim et al.(2) believe that perishable food has a relatively short shelf life and begins to spoil shortly after production. The perishability feature dramatically affects the production strategy and planning process by adding new requirements. Since food products are classified as highly perishable goods, the freshness of these products is significantly influenced by the time and temperature of the environment during the delivery. Therefore, keeping them fresh during the delivery time is very difficult. In this regard, one of the most effective solutions is to apply the vehicle routing problem (VRP) to deliver food products to customers in the shortest possible time and with the minimum cost. Many researchers in mathematical models have considered different types of VRP problems since Dantzing and Ramser (3) first introduced them. Vehicle routing problems are optimization problems that occur in the final phase of the supply chain. The approach to achieving the best solution for these problems is essential for protecting natural resources, saving energy, and increasing the efficiency of logistics activities. In these problems, a feasible solution includes a series of routes so that each customer is visited only once and exactly by one vehicle. The simplest and simultaneously most suitable type of VRP problem is the capacitated vehicle routing problem (CVRP). In these problems, the vehicles are identical and located in a single central depot. Each vehicle can be assigned at most one route, and the quantity of goods delivered on each route does not exceed the vehicle capacity. In these problems, the objective function is usually used to minimize costs or maximize the total profit(4). Another debatable issue that raises the market's competitiveness level is the attention to significant changes in the business environment, such as customer demand, which has attracted researchers to design a robust supply chain (5).

In this research, we investigate the production and distribution planning, scheduling, and routing in the yogurt product supply chain in the frame of an integrated mathematical Mixed-integer non-linear programming model to account for two arguably important factors: perishability and uncertainty. So, we can give the decision-makers a proper perspective to help them take steps in planning confidently. This problem with considering capacity constraints in the production network and transporting products to the distribution centers has characteristics that provide the conditions for its application in the dairy industry. Uncertain demand during the planning horizon, production constraints, the perishability feature of the product, sequence dependence, quality control process, packaging, and distribution, are some of these characteristics. On the other hand, delivering products as freshly as possible is vital to customer satisfaction. A meta-heuristic approach based on a genetic algorithm is also introduced to solve the problem. Its efficiency is investigated in a wide range of sample problems and a real case study. The results show the effectiveness of the proposed solution in solving different problems.

The subsequent parts of the research are organized as follows. In section 2, the relevant literature is reviewed. A brief description of the yogurt production process is considered in section 3. Afterward, section 4 describes the mathematical formulation of the integrated planning, scheduling, and routing model. Section 5 aims to define the robust fuzzy programming approach. The proposed solution approach is described in section 6, and an illustrative case study is discussed with computational results in section 7. Finally, some concluding remarks are drawn, and future research lines are given in section 8.

# 2. RESEARCH BACKGROUND

This section will review the most relevant and closest literature to this research and introduce outstanding articles. In today's world, the food supply chain (FSC) plays an essential role in our daily lives. FSC management has attracted much attention from universities, and many efforts have been made to raise the standard of living. However, there are many issues, such as high supply chain costs, food waste, and environmental issues in the food industry. Hence, researchers have challenged various aspects of this industry. Ding et al. (6) examined factors that appoint the competitive advantage of dairy supply chains using an evidence-based approach with the help of evidence from the Chinese dairy industry. Their main focus was on the quality assurance of dairy products, which they believed was one of the influential fundamental factors. Jouzdani and Govindan (7) studied a multi-objective mathematical programming model to optimize the cost, energy consumption, and traffic congestion associated with the sustainability of food supply chain operations. They considered social and environmental impacts along with their well-known economic aspects. The results indicated that emphasizing the economic aspect, the chain's environmental impact may increase by 120% for highly perishable products. The social impact may rise by 51% for the highly congested road networks. However, a 15% economic compromise can improve the sustainability of the supply chain network design by 150%. Tostivint et al. (8) believed that food waste during the supply chain is too much to be ignored. They mentioned that the United Nations Environment Program and the World Resources Institute launched the Global Food Loss & Waste (FLW) Protocol in October 2013. Their goal was to foster an "FLW Standard" for accounting and reporting amounts of food waste within the food supply chain. Ledo et al. (9) stated that the demand for milk and milk products is increasing at a high rate, which may result in rapidly evolving dairy chains. Because they were not satisfied with the standard of milk quality and safety, their study aimed to get insight into possible causes of persisting poor milk safety and hygiene practices. The Tanzanian dairy chain was also taken as a case study.

Recent reviews covering the issues in supply chain planning highlight the importance of corruption considerations. Yavari and Geraeli (5) designed a green closed-loop supply chain network for perishable products. Some parameters, such as demand, return rate, and the quality of returned products, are considered uncertain. They developed a mixed-integer Linear programming model to achieve the objective of minimizing total cost and environmental pollutants. They also adopted an innovative MILP robust model to address the uncertainty. Aras and Bilge (10) introduced the problem of designing a Robust Supply Chain in the snacks market for a company in Turkey. According to the prosperity of this market, a new strategy has been formulated to raise a new plant. At the decision-making level, they address the problem of timing and location of the new plant, its initial capacity, the allocation of demand, and the quantity of raw material delivered from the suppliers in each period. The object of their problem is to minimize total cost, and they develop the problem using a mixed-integer Linear programming model, which is solved by the commercial solver CPLEX. Banasik et al. (11) introduced a stochastic programming model with multi-objective and two-stage to analyze the economic and environmental effects and its role in uncertainty in agri-food supply chains. A mushroom supply chain in the Netherlands is addressed as a case study. In the end, they compared the optimal production planning decisions made by a two-stage stochastic programming model with the results of an equivalent deterministic model. Rahbari et al. (12) investigated a bi-objective MILP model to consider the vehicle routing and scheduling problem with cross-docking for perishable products. Assuming that the travel time of the vehicles and the shelf life of the products is uncertain, two robust models are introduced. Rafie-Majd et al. (13) studied a three-echelon supply chain modeled as the integrated inventory location routing problem (ILRP). They considered the limited time horizon in which perishable products are delivered to the customer. The retailers have a stochastic demand and follow a normal distribution with a particular mean and standard deviation. A Lagrangian Relaxation Method is adopted to solve the problem. Mohammadi et al. (14) developed a novel mechanism in fresh produce supply chain management aiming to decrease food waste during the supply chain and improve the profits of independent food supply chain members. Banasik et al. (11) adopted a multi-objective two-stage stochastic programming in the agri-food supply chain. To consider uncertainty, they analyzed the economic and environmental issues. The proposed model is applied in a mushroom supply chain in the Netherlands. Vahdani et al. (15) addressed the problem of production, inventory, and routing for deteriorating products. They divided their problem into two phases production and delivery. They considered a mathematical programming approach to maximize profits, which is the objective of their model. Two heuristic and meta-heuristic algorithms are proposed to solve the problem. Jouzdani et al. (16) studied multiple products, multiple transportation modes, the monetary value of time, and uncertainty in transportation costs, demand, and supply. Hooshangi-Tabrizi et al. (17) presented an inventory management problem with perishability and demand uncertainty where the goal is to minimize the sum of ordering, purchasing, holding, shortage, wastage, and modification costs. They formulated the problem as a novel two-stage robust integer optimization model and developed an exact column-and-row generation algorithm to solve it. Alinezhad et al. (18) introduced a sustainable closed-loop supply chain network configured under uncertain conditions. The proposed network is a multi-product, multi-period problem in which the objectives are to maximize the supply chain profit and customer satisfaction at the same time. Moreover, the carbon footprint is included in the first objective function regarding cost (tax) to affect the total profit and treat the environmental aspect. Yantong LI et al. (4) formulated the production inventory routing planning with an integrated mixed integer linear programming (MILP) model where the food quality level is explicitly traced throughout the supply chain, and the objective of their proposed model is to maximize the total profit. According to Amorim et al. (19) and (20), shelf life can be considered a loss or benefit function to consider the economic value of fresh products in objective functions. However, some other researchers have also considered shelf life as a constraint (21), (1), (22), and (23). We can also ignore the Shelf-life effect on the objective function and consider it only as a constraint (2).

Much research has also been conducted to find an efficient routing policy for perishable goods. Meidute-Kavaliauskiene et al. (24) addressed a VRP multi-period and multi-echelon perishable supply chain problem concerning procurement time, cycle cost, and customer satisfaction. This study presented a new form accounting for environmental considerations, cost, procurement time, and customer satisfaction, such that the total costs, delivery time, and the emission of pollutants in the network were minimized while customer satisfaction was maximized. Suryawanshi et al. (25) contributed to the scant literature on the quantitative modeling of food SC. The research work integrated a stochastic nature of SC simultaneously coupled with the effect of disruption, transport losses, and product perishability. It incorporated proactive planning strategies to minimize the disruption impact and the concept of incremental quantity discounts on lot sizes at a destination node. Majidi et al. (26) presented an integrated multi-objective sustainable pricing-production-workforce-routing problem for perishable products. Total profit, workforce planning, and vehicle fuel consumption are considered objective functions due to the importance of operational performance and social and environmental concerns. The application of the proposed approach is investigated using actual case data from a dairy product supply chain. In addition, it concurrently optimizes the selling price and protects the environment from the negative impacts of greenhouse gas emissions (GHGs). Hasani et al. (27)designed a closed-loop supply chain for corruptible goods. This paper considers multiple periods, multiple products, and a multi-level supply chain with uncertain demand. K. Govindan et al. (28) offered a multi-vehicle routing problem for a two-level supply chain for perishable food products. They developed a multi-objective optimization model for the perishable food supply chain network, which involves determining the number and location of the factories and optimizing the number of products shipped in the shortest path at each level. This research uses a heterogeneous type of vehicle for delivering products. In addition, some research has focused on the distribution of frozen and cold foods, commonly referred to as the cold supply chain. Duk Song and Dae Ko (29) focused on a real-life routing problem of a multi-product model for the daily delivery of perishable food products. To maximize the total satisfaction of customers, which is the result of the freshness of the delivered product, they simultaneously used both refrigerated and conventional vehicles. Bilgen and Celebi (22) presented a MILP model addressing the production scheduling and distribution planning problem in a yogurt production line of multi-product dairy plants. They consider yogurt production with perishability and sequence-dependency issues by focusing on the packaging stage operating with parallel units sharing common resources also, Bilgen and Celebi (22), in their work, considered a production and distribution problem for set-type yogurt production. The model has been extended by considering timing and capacity constraints concerning the incubation operation of set-type yogurt. Mousavi et al. (30) recently presented a new production routing model for perishable products with uncertain demand. The aim is to minimize the costs of production, inventory, routing, wasted products, and penalties for non-fresh products. The model is more applicable for perishable products with a limited, discrete shelf life and a high freshness value. Computational experiments show that the proposed mathematical model can significantly reduce wasted products, particularly when consumer buying patterns change due to various occurrences, such as a pandemic. Pratap et al. (31) presented a framework to create an integrated production-inventory-routing problem for perishable food products that considered capacity, time windows, and carbon emissions reduction. As the inventory routing problems are NP-hard, the proposed problem was solved using two nature-inspired algorithms whose computational results indicated mitigation in the stochastic nature of problems in uncertain circumstances. Table 1 depicts an overview of the related literature. Hashemi-Amiri et al.(32) To solve an integrated supplier selection, production planning, and sourcing problem, a bi-objective optimization
model for securing shipments through PFSC, which has three levels of perishable food supply chain is proposed. This network includes multiple products. By focusing on optimal network costs, the proposed model seeks to decrease supply and demand uncertainty and enhance distribution decision-making**.**

**TABLE 1.** Overview of the related literature

| Research | Production | Scheduling | Storage | Distribution | Multi Product | Working Time | Perishability | Uncertainty | VRP | Case Study |
|---|---|---|---|---|---|---|---|---|---|---|
| Hasani et al. (2012) | | | ✓ | ✓ | ✓ | | | | ✓ | |
| Bilgen and celebi (2013) | ✓ | ✓ | ✓ | ✓ | ✓ | ✓ | ✓ | | | |
| Sel et al. (2015) | ✓ | ✓ | ✓ | ✓ | ✓ | ✓ | ✓ | | | ✓ |
| Yantong LI et al. (2016) | ✓ | | ✓ | ✓ | | | ✓ | | ✓ | |
| Vahdani et al. (2017) | ✓ | ✓ | ✓ | ✓ | | | ✓ | | ✓ | |
| Yavari and Geraeli (2019) | ✓ | | ✓ | | ✓ | | ✓ | ✓ | | ✓ |
| Jouzdani et al. (2020) | ✓ | ✓ | | ✓ | ✓ | | ✓ | ✓ | | ✓ |
| Jouzdani and Govindan (2021) | ✓ | ✓ | ✓ | ✓ | ✓ | ✓ | ✓ | ✓ | | ✓ |
| Hooshangi-Tabrizi (2022) | | | ✓ | ✓ | ✓ | | ✓ | ✓ | | ✓ |
| Meidute-Kavaliauskiene et al. (2022) | | ✓ | ✓ | ✓ | ✓ | | ✓ | | ✓ | ✓ |
| Alinezhad et al. (2022) | ✓ | | ✓ | ✓ | ✓ | | ✓ | ✓ | | ✓ |
| Majidi et al. (2022) | ✓ | | ✓ | ✓ | ✓ | ✓ | ✓ | | | ✓ |
| Pratap et al. (2022) | ✓ | ✓ | ✓ | ✓ | | | ✓ | ✓ | ✓ | |
| Mousavi et al. (2022) | ✓ | | ✓ | | | | ✓ | ✓ | | |
| Hashemi-Amiri et al.(2023) | ✓ | ✓ | | | ✓ | | | ✓ | ✓ | ✓ |
| This study(2024) | ✓ | ✓ | ✓ | ✓ | ✓ | ✓ | ✓ | ✓ | ✓ | ✓ |

## 3. DESCRIPTION OF THE YOGURT PRODUCTION PROCESS

Yogurt production is used to typify the production of all fermented products. The milk used is pre-treated according to the product to be manufactured and undergoes a standardization process to ensure a specific fat content. Regardless of the yogurt type, be it set yogurt, stirred yogurt, drinking yogurt, or another fermented product such as sour cream, kefir, buttermilk, or sour milk, the core process is the same. The milk is brought to the required fat content for the desired end product. After homogenization and pasteurization, the bacterial culture specific to the product is added to the milk and incubated. Once products reach the optimum pH value, they are cooled, a fruit mixture is added if required, and then filling takes place. This process is usually under aseptic conditions to avoid recontamination.

Among the several types of yogurts, the two main types are set and stirred yogurt. *Stirred yogurt* is a common type that is usually combined with several elements such as fruit, nuts, and other ingredients additions. Hence, we are confronted with many products in the production of stirred yogurt. For this reason, producing this type of yogurt is much more challenging than other types. The yogurt production process starts with the collection of milk and continues with pasteurization, standardization, homogenization, culture addition, fermentation/incubation, packaging, and cold storage and distribution operations (23). According to Chandan and Shahani (33), the success of the yogurt industry can be attributed to its health-related glamour, the increase of low-fat products, the wide variety of tastes and textures, intense marketing activities, and its relatively low costs. Figure 1. Shows a diagram of yogurt production flow.

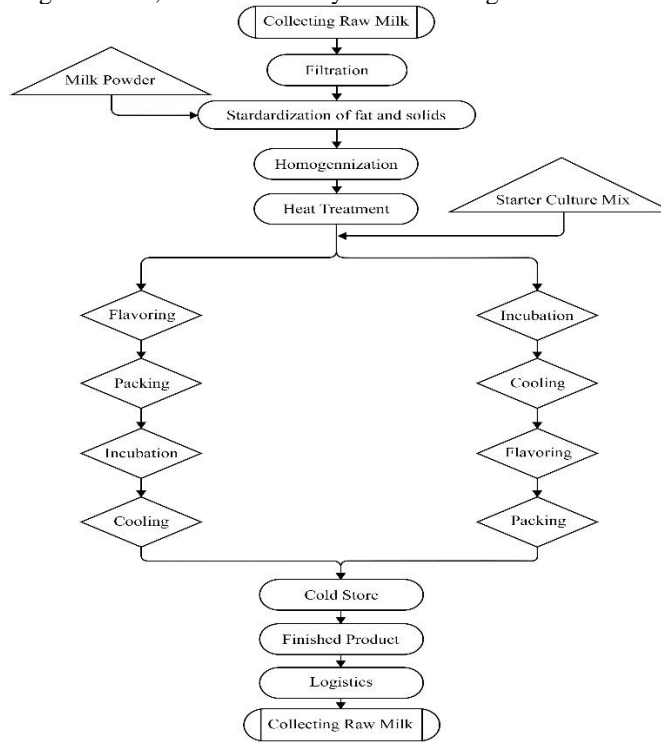

**Figure 1.** Yogurt production line

## 4. MATERIALS AND METHODS

The proposed model in this research is inspired by a real-case dairy industry (Kalleh, Tehran) producing yogurt products. A Mixed-integer non-linear programming (MINLP) framework is developed based on the definition of families of products. Our integrated MINLP model is an extension of the formulation proposed by Bilgen and Celebi (22) and the problem of Simultaneous production and logistics operations planning in semi-continuous food industries proposed by Georgios M. Kopanos (34). The problem involves Tactical and operational decisions and vehicle routing.

The key decisions for each planning and scheduling period are as follows:
1. The optimal assignment of families and products to processing lines in each production period.
2. The family sequencing on each line in each period.
3. The inventory level of each product in each period.
4. The produced amount of each product in each line at each period.
5. The quantity of transferred product to each distribution center and the corresponding unmet demand quantity.
6. Processing time of each family on each line, and overtime.
7. Maximum completion time of each family and each line.
8. Allocating each vehicle to each route to deliver the final product to each distribution center

Finally, the objective function is to fully satisfy customer demand with the lowest total cost, including production costs, set-ups, changeovers, inventory, overtime, unmet demand, and transportation.

## 4.1. Problem-Assumption

In this section, we officially address the problem of production planning and logistics operations in a semi-continuous food process industry with multiple products. The essential features of the problem under consideration are summarized as follows. Some of these characteristics are adapted from the assumptions of two papers, Bilgen and Celebi (22) and Georgios M. Kopanos (34).

The supply network includes a single plant that delivers the final products to several distribution centers.

Demand is an uncertain parameter.

The plant includes a set of products with many similar features, which allows us to group products with similar characteristics into a product family. Therefore, the production planning problem could be partially focused on families rather than on each product separately. Specifically, sequencing and timing decisions are taken for families rather than each product. The demand for each product in each period is collected from the distribution centers. The demands have specific due dates, and backlogging is not allowed. Unmet demand cannot be transferred to the following periods. The unmet demand is discarded at costs.

A known planning horizon is divided into a set of periods.

A set of products with uncertain demand in each period, inventory costs, production rates, variable production costs, and minimum required time for storing and cooling processed products are available. Further, the proposed modeling approach allows products from the same family to have different processing rates, operating costs, setup times, inventory, and transportation costs.

The planning and scheduling are performed on a short-term horizon. For each product, inventory balances are updated daily according to the production output from the plant. There is a specific inventory capacity for final products, which are stored in pallets. Storage costs depend on product types.

The plant includes several identical production lines with available processing time in each period. Each line can produce all products and has different operating costs. A product cannot be processed on multiple lines simultaneously, and a line cannot process more than one job simultaneously. Lines are always available. The production is limited to the minimum and maximum lot sizes.

After finishing the production process, they place packed products in cooling storage containers for about 1–3 days to achieve ultimate stability and preserve a high quality for the final product. Quality control is also addressed at this phase. The products are not allowed to be delivered before they complete the required time for cooling storage and the quality control processes.

Freshness plays an essential part in the competition. Freshness can be measured with shelf life which is a small duration. The minimum shelf life required by the customers is defined as a critical rate which is a fraction of the maximum shelf life.

The available working hours for the lines are defined according to working days. Under conditions of heavy demand, overtime production is allowed every working day.

A set of batch recipes (fermentation recipes) with minimum preparation time, preparation cost, and minimum and maximum production capacity are available.

Changeover operation is required on each processing unit whenever the production is changed between two families. For each changeover operation, time and cost are considered separately.

The plant includes a set of vehicles that transfer the final product to the distribution centers under the following circumstances:

A min and max capacity characterize each vehicle.

Each car is assigned at most one route.

Each customer is served by exactly one vehicle.

Each route begins and ends at the depot in the factory.

At most, the total vehicle load at any route point is the vehicle capacity.

## 4.2. Notation

The following notation is applied to formulate the mathematical problem:

**TABLE 2.** Indicators parameters and variables of the proposed model

| Notation | Description |
| --- | --- |
| Sets | |
| $a,b,c$ | Distribution Centers (DCs) ($a,b,c = 1...A$) |
| $f$ | Product family ($f = 1...F$) |
| $j$ | Production line ($j = 1...J$) |
| $l$ | Vehicle ($l = 1...L$) |
| $p$ | Product ($p = 1...P$) |
| $d$ | Demand day ($d = 1...D$) |
| $r$ | Batch recipe ($r = 1...R$) |

| | $i$ | Production day ($i = 1...I$) |
|---|---|---|
| Sub-Sets | | |
| | $F_j$ | families that can be processed on unit $j$ |
| | $J_f$ | available units $j$ to process family $f$ |
| | $J_p$ | units $j$ that can process product $p$ |
| | $P_f$ | products $p$ that belong to the same family |
| | $R_f$ | batch recipes $r$ for family $f$ |
| | $R_j$ | batch recipes $r$ that can be processed on unit $j$ |
| | $F_r$ | family $f$ that has the same batch recipe $r$ |
| Parameter | | |
| | $\Pr et_r$ | minimum time for preparing batch recipe (min) |
| | $Rtime_i$ | regular working time on day $i$ (min) |
| | $Maxtime_i$ | Maximum available time on day $i$ (min) |
| | $ShelfLife_p$ | shelf life of product $p$ (d) |
| | $CrRate_p$ | minimum shelf life requirement of customer for product $p$ (% of shelf life) |
| | $FCT_l$ | fix cost for contracting transportation truck $l$ |
| | $Var\,Cost_p$ | the variable production cost of product $p$ |
| | $Over\,Cost_i$ | overtime cost (min) |
| | $MaxTC_l$ | maximum capacity of transportation truck $l$ (kg) |
| | $MinTC_l$ | minimum capacity of transportation truck $l$ (kg) |
| | $Pcapacity_i$ | production capacity of the plant (kg) |
| | $Pallet_p$ | the factor for converting product to the storage unit |
| | $dailyop_{ji}$ | daily opening setup time for every unit $j$ on the day $I$ (min) |
| | $dailysh_{ji}$ | daily shutdown time for every unit $j$ on day $i$ (min) |
| | $W_{ji}$ | cost for producing batch recipe $r$ on the day $i$ |
| | $Cht_{fej}$ | changeover time between family $f$ and family $e$ in unit $j$ on day $i$ (min) |
| | $Setup_{jp}$ | setup time for product $p$ in unit $j$ (min) |
| | $Relt_{ri}$ | release time for batch recipe $r$ on day $i$ (min) |
| | $O_{ji}$ | additional unit preparation time for processing unit $j$ on day $i$ (min) |
| | $Bpc_{ri}$ | cost for producing batch recipe $r$ on the day $i$ |
| | $IC_{pi}$ | inventory cost of product $p$ (D/l per day) |
| | $VTC_{abl}$ | cost for transferring products from distribution center $a$ to distribution center $b$ by truck $l$ (Km/h) |
| | $Chc_{feji}$ | changeover cost between family $f$ and family $e$ in unit $j$ on day $i$ |
| | $Line\,Cost_{ji}$ | operating cost of line $j$ on the day $i$ (per day) |
| | $MaxLots_p$ | maximum production lots of product $p$ (kg) |

| | |
|---|---|
| $MinLots_p$ | minimum production lots of product $p$ (kg) |
| $Demand_{dpa}$ | the demand of distribution center $a$ for product $p$ on demand day $d$ (kg) |
| $Prate_{jp}$ | production rate of product $p$ on line $j$ |
| $UnmdCost_{ap}$ | the cost of unmet demand of product $p$ for DC $a$ |
| $FCost_{ji}$ | fixed cost for utilizing unit $j$ on the day $i$ |
| $MuMax_{ri}$ | maximum production capacity of batch recipe $r$ on day $i$ (kg) |
| $MuMin_{ri}$ | minimum production capacity of batch recipe $r$ on day $i$ (kg) |
| $CST$ | Minimum required time to cool the finished product before delivering (min) |
| $QCTime$ | required time for quality control (min) |
| $StCapacity$ | storage capacity of the plant |
| $M$ | Big number |
| $\alpha$ | minimum satisfaction level of the flexible constraint |
| $\gamma$ | the penalty value of the possible violation of the soft constraint |

Decision variables

| | |
|---|---|
| $UD_{apld}$ | amount of product $p$ delivered to distribution center $a$, by vehicle $l$ on demand day $d$ (kg) |
| $UV_{abpld}$ | amount of product $p$ transferred and delivered to distribution center $b$, by vehicle $l$, after visiting distribution center $a$ on-demand day $d$ (kg) |
| $UB_{abld}$ | the total amount of products transferred to distribution center $b$, by vehicle $l$, after visiting distribution center $a$ on-demand day $d$ (kg) |
| $Q_{pji}$ | produced amount of product $p$ in unit $j$ on day $i$ (kg) |
| $QB_{pi}$ | total produced amount of product $p$ on day $i$ (kg) |
| $UnmD_{adp}$ | unmet demand of product $j$ for DC $a$ on-demand day $d$ (kg) |
| $II_{pi}$ | inventory of product $p$ on day $i$ (kg) |
| $PT_{fji}$ | processing time for family $f$ in unit $j$ on day $i$ (min) |
| $CT_{fji}$ | completion time for family $f$ in unit $j$ on the day $i$ (min) |
| $Cmax\text{-}Family_{fi}$ | maximum completion time for family $f$ on day $i$ (min) |
| $Cmax\text{-}Line_{ji}$ | maximum completion time of unit $j$ on day $i$ (min) |
| $OverTime_i$ | required overtime on day $i$ (min) |
| $V_{ji}$ | equal to 1, if unit $j$ is used on day $i$ |
| $G_{ri}$ | equal to 1, if batch recipe $r$ is produced on day $i$ |
| $X_{feji}$ | equal to 1, if family $e$ is processed exactly after family $f$, when both are assigned to the same unit $j$ on day $i$ |
| $Y_{fji}$ | equal to 1, if family $f$ is assigned to unit $j$ on the day $i$ |
| $YB_{pji}$ | equal to 1, if product $p$ is assigned to unit $j$ on day $i$ |
| $ZV_{abld}$ | equal to 1, if transportation truck $l$ transfers product from distribution center $a$ to distribution center $b$ on demand day $d$ |
| $Z$ | The value of the total objective function |

### 4.3. Problem Formulation

According to the above notation, the MINLP problem is formulated as follows:

$$Min\ Z = \sum_{ipj} Q_{pji} * VarCost_p + \sum_{pi} II_{pi} * IC_{pi} + \qquad (1.1)$$
$$\sum_{ri} BPC_{ri} * G_{ri} + \sum_{feji} X_{feji} * Chc_{feji} + \qquad (1.3)$$
$$\sum_{i} OverTime_i * OvertCost_i + \qquad (1.5)$$
$$\sum_{adp} UnmD_{adp} * UnmDCost_{ap} + \qquad (1.7)$$
$$\sum_{ij}(FCost_{ji} + LineCost_{ji}) * V_{ji} +$$
$$\sum_{ald,\,a\neq 1} FCT_l * ZV_{1ald} + \sum_{abld,\,a\neq b}(UB_{abld} * VTC_{abl}) * ZV_{abld}$$

(1.1) (1.2) (1.3) (1.4) (1.5) (1.6) (1.7) (1.8)

In Eq. (1), the non-linear objective function aims at minimizing total cost. The total cost corresponds to variable production cost(Eq.(1.1)), inventory cost (Eq.(1.2)), batch recipes preparation cost(Eq.(1.3)), family changeover cost(Eq.(1.4)), overtime cost(Eq.(1.5)), unmet demand cost(Eq.(1.6)), line utilization cost(Eq.(1.7)), and transportation cost (Eq.(1.8)).

$$II_{pi} = QB_{pi} - \sum_{a,a\neq 1}\sum_{l} UV_{1apli} \qquad \forall i=1, p, d=i \qquad (2)$$

$$II_{pi} = II_{p(i-1)} + QB_{pi} - \sum_{a,a\neq 1}\sum_{l} UV_{1apli} \qquad \forall i>1, p, d=i \qquad (3)$$

Eq. (2) shows the inventory level only for the first day. The inventory of product $p$ at the end of the first day equals the difference between the quantity of product $p$ produced on the first day and the delivered quantity of product p produced to the DCs. Eq. (3) refers to the inventory level of product $p$ at the end of the day $i$. The inventory level is computed by adding the production quantity of product $p$ on the day $i$ to the difference between the inventory of the previous day and the delivered quantity of product $p$ to the DCs on the same day.

$$Rtime_i + Overtime_i \leq Maxtime_i \qquad \forall i \qquad (4)$$
$$Cmaxfamily_{fi} - Rtime_i \leq Overtime_i * M \qquad \forall i, f \qquad (5)$$
$$Cmaxfamily_{fi} \leq Overtime_i + Rtime_i \qquad \forall i, f \qquad (6)$$

Eq. (4) limits the total working time with the maximum available time on the day $i$. Eq. (5) limits the maximum completion time on day $i$ with the regular working time if overtime is unnecessary. Eq. (6) shows that the maximum completion time on day $i$ is limited, with maximum working time consisting of regular hours and required overtime.

$$cmaxfamily\ fi = \sum_{j\in J} CT_{fji} \qquad \forall i, f \qquad (7)$$
$$cmaxfamily\ f_i \geq \sum_{j\in J} PT_{fji} \qquad \forall i, f \qquad (8)$$

Eq. (7) and (8) present that the maximum Completion time of a family on the day $i$ must be equal to the summation of completion time of that family on all production lines ($j\epsilon J_f$) and equal to or greater than the summation of the processing time of that family on all production lines ($j\epsilon J_f$), respectively
.

$$MuMin_{ri} * G_{ri} \leq \sum_{p \in P_r} \sum_{j \in J_p} Q_{pji} \qquad \forall i, r \qquad (9)$$

$$\sum_{p \in P_r} \sum_{j \in J_p} Q_{pji} \leq MuMax_{ri} * G_{ri} \qquad \forall i, r \qquad (10)$$

Batch recipe stage constraints must be included in the mathematical model to ensure a feasible production plan, according to Eqs. (9) and (10) the cumulative produced quantity of products $p \in P_r$ should be greater than the minimum production capacity of the batch recipe and lower than its maximum production capacity.

$$\sum_{f \in F_j} \sum_{\substack{e \in F_i \\ e \neq f}} X_{feji} + V_{ji} = \sum_{f \in F_j} Y_{fji} \qquad \forall j, i \qquad (11)$$

$$V_{ji} \geq Y_{fji} \qquad \forall f, i, j \in J_f \qquad (12)$$

$$YB_{pji} * Minlots_p \leq Q_{pji} \qquad \forall i, p, j \in J_f \qquad (13)$$

$$Q_{pji} \leq YB_{pji} * Maxlots_p \qquad \forall i, p, j \in J_f \qquad (14)$$

Eq. (11) defines that the total number of active sequencing binary variables plus the line utilization binary variable should equal the total number of active binary variables according to assigning the family products to available production *line j*.
Eq. (12) ensures that the processing line *j* is used at day *i* if at least one family *f* is assigned to it over this period. Eq. (13) and (14) define minimum and maximum production lot sizes.

$$Y_{fji} \geq YB_{pji} \qquad \forall f, i, j \in J_f, p, \in P_f \qquad (15)$$

$$Y_{fji} \leq \sum_{p \in p_F} YB_{pji} \qquad \forall f, i, j \in J_f \qquad (16)$$

$$\sum_{\substack{e \in F_j \\ e \neq f}} X_{feji} \leq Y_{fji} \qquad \forall f, i, j \in J_f \qquad (17)$$

$$\sum_{\substack{e \in F_j \\ e \neq f}} X_{efji} \leq Y_{fji} \qquad \forall f, i, j \in J_f \qquad (18)$$

$$QB_{pi} = \sum_{j \in J_p} Q_{pji} \qquad \forall p, i \qquad (19)$$

Eq. (15) A family is assigned to a processing line j on day i if at least one product that belongs to this family is processed in this line during the same day. So, Eq. (16) enforces the binary variables $Y_{fji}$ to zero when no products are processed in line j during day *i*. Eq. (17) and (18) state that if a family f is allocated to processing line j on day *i*, then at most, one family is processed before and after it, respectively. The total quantity of product *p* produced on the day i is calculated as Eq. (19).

$$PT_{fji} = \sum_{p \in P_f} \left( \frac{Q_{pji}}{Prate_{pj}} + Setup_{jp} * YB_{pji} \right) \qquad \forall f, i, j \in J_f \qquad (20)$$

Because products that belong to the same family do not require changeover operations among them, sequencing and timing constraints should be only considered on families. To define sequencing and timing decisions for families, the definition of family processing time is introduced as Eq. (20).

$$Mintc_l * ZV_{a1ld} \leq UB_{a1ld} \qquad \forall a \neq 1, l, d \qquad (21)$$

$$UB_{a1ld} \leq Maxtc_l * ZV_{a1ld} \qquad \forall a \neq 1, l, d \qquad (22)$$

$$\sum_{a \neq 1} ZV_{a1ld} = 1 \qquad \forall l, d \qquad (23)$$

$$\sum_{b \neq a} Z_{abld} = \sum_{b \neq a} Z_{bald} \qquad \forall a, l, d \qquad (24)$$

$$\sum_a \sum_l ZV_{abld} = 1 \qquad \forall d, b \neq a, l \qquad (25)$$

$$\sum_{a \neq 1} ZV_{a1ld} = 1 \qquad \forall d, l \qquad (26)$$

$$UB_{abld} = \sum_p UV_{abpld} \qquad \forall d, l, a, b, a \neq b \qquad (27)$$

$$UD_{apld} \leq UV_{abpld} \qquad \forall a \neq b \neq 1, p, l, d \qquad (28)$$

$$UB_{a1ld} \leq M * ZV_{a1ld} \qquad \forall a \neq 1, l, d \qquad (29)$$

$$UB_{abld} \leq M * ZV_{abld} \qquad \forall a \neq b \neq 1, l, d \qquad (30)$$

$$o_a - o_b + (M)ZV_{abld} \leq M - 1 \qquad \forall a, b \in N \setminus \{1\} \text{ and } \forall a \neq b \text{ and } \forall l, b \qquad (31)$$

According to Eq. (21) and (22), each vehicle has a specific minimum and maximum capacity that needs to be considered. Eq. (23) means that if you use vehicle l, travel must be started from the depot. Eq. (24) represents the flow balance (if the vehicle enters a node, it exits from that node). Eq. (25) implies that any distribution center with demand should be covered. Eq. (26) ensures that vehicles brought to the distribution centers from the depot must be returned to the depot. The total amount of transported products is calculated using Eq. (27). Eq. (28) ensures that the number of products delivered to a distribution center will surely be less than or equal to the inventory level of that product in the truck. According to Eq. (29) and (30), ensure that, from one center to another, the products are transported by a vehicle only if a trip between those two centers is made. Eq. (31) is the sub-tour elimination constraint.

$$CT_{fji} - PT_{fji} \geq (dailyop_{ji} + Pret_r + Max[O_{ji}, Relt_{ri}]) * Y_{fji} + \sum_{\substack{e \in F_j \\ e \neq f}} Cht_{efj} * X_{efji} \qquad \forall f, i, j \in J_f, r \in R_f \qquad (32)$$

$$CT_{fji} \leq (W_{ji} - dailysh_{ji}) * Y_{fji} \qquad \forall f, i, j \in J_f \qquad (33)$$

$$CT_{fji} + Cht_{fej} \leq CT_{eji} - PT_{eji} + M * (1 - X_{feji}) \qquad \forall f, i, j \in J_f, f \neq e, e \qquad (34)$$

Eq. (32) and (33) refer to the beginning and finishing time of each family of products which is adopted from Georgios M. Kopanos [39]. Eq. (34) states that the starting time of family e, which directly follows another family f on a processing line j, should be greater than the completion time of family f, plus the necessary changeover time between those families.

$$G_{ri} \geq \sum_{j \in J_f} Y_{fji} \qquad \forall r, i, f \in F_r \qquad (35)$$

$$CT_{fji} \leq Cmaxline_{ji} \qquad \forall r, i, f \in F_j \qquad (36)$$

Eq. (35) Shows that a batch recipe is produced on the day $i$ if at least one family is processed online $j$ on same day. Eq. (36) defines that the completion time for the families processed on line $j$ cannot exceed the maximum available time of line $j$.

$$\sum_p II_{pi} * Pallet_p \leq Stcapacity \qquad \forall i \qquad (37)$$

$$\sum_p QB_{pi} \leq Pcapacity_i \qquad \forall i \qquad (38)$$

Eq. (37) converts the inventory of product *p* on day *i* into unit storage (pallets) by multiplying it with the corresponding factor and limiting the stored inventory level. The plant has a production capacity that, according to Eq. (38), the total amount of products produced cannot exceed this capacity.

$$\sum_{f \in F_j} PT_{fji} + \sum_{f \in F_j} \sum_{\substack{e \in F_j \\ e \neq f}} Cht_{fej} * X_{feji} \leq \left(W_{ji} - dailysh_{ji} - dailyop_{ji} - Min_{r \in R_j}[Pret_r]\right) * V_{ji} \qquad \forall i, j \qquad (39)$$

$$PT_{fji} + \sum_{\substack{e \in F_j \\ e \neq f}} Cht_{fej} * X_{feji} \leq \left(W_{ji} - dailysh_{ji} - dailyop_{ji} - Pret_r\right) * Y_{fji} \qquad \forall f, i, j \in J_f, r \in R_f \qquad (40)$$

To reduce the computational effort, Eq. (39) can impose an upper bound on the total processing time for each processing line at each day *i*. Similarly, in Eq. (40), an upper bound on the family processing time can be defined.

$$Demand_{dpa} \leq \sum_l \sum_b UD_{bapld} * ZV_{bald} + UnmD_{adp} \qquad \forall a \neq 1, d, p, b \neq a \qquad (41)$$

$$\sum_{a \neq 1} \sum_l \sum_p UV_{1apld} \leq \sum_p II_{pi} \qquad \begin{array}{l} \forall i, d, p, i = d - S \\ d - i < (1 - CrRte_p) * Shelflife_p \end{array} \qquad (42)$$

Eq. (41) states that the demand of each center for each product is less than the total amount of that product transported to all distribution centers plus the amount of unmet demand for that product. Following Eq. (42), the quantity of product shipped to distribution centers on day *d* should be less than the inventory level of that product on day *i = d - S* due to the necessity for quality control and a cooling period.

$$\begin{array}{l} UB_{abld}, UD_{apld}, UV_{abpld}, Q_{pji}, QB_{pi} \\ UnmD_{adp}, II_{pi}, PT_{fji}, CT_{fji}, Cmaxfamily_{fi} \\ Cmaxline_{ji}, Overtime_i \geq 0 \\ V_{ji}, G_{ri}, X_{feji}, Y_{fji}, YB_{pji}, ZV_{abld} \in \{0, 1\} \end{array} \qquad \forall i, d, p, a, f, j, i \qquad (43)$$

Finally, Eq. (43) shows the non-negativity and integrality features of the variables.

### 4.4. The Linearized Mixed Integer Programming Formulation

The transformed mixed integer linear programming formulation is referred to as the proposed model and is given below.

Minimize Eq.(1.1)+Eq.(1.2)+Eq.(1.3)+Eq.(1.4)

$$+ Eq.(1.5) + Eq.(1.6) + Eq.(1.7)$$

$$+ \sum_{ald, a \neq 1} FCT_l * ZV_{1ald} + \sum_{abld, a \neq b} (UB_{abld} * ZVTC_{abl})$$

Subject to Equations (2)-(44).
The linearization procedures implemented within the proposed model can be found in the proposition below.

**Proposition:** Now we consider a change of variables, let us introduce a new continuous variable $ZVTC_{abld}$ which equals multiplying both variables $VTC_{abl}$ and $ZV_{abld}$ as given in the below equation.

$$\sum_{abld, a \neq b} (UB_{abld} * VTC_{abl}) * ZV_{abld}$$

This nonlinear component can be linearized under the following set of constraints:

$$ZVTC_{abld} \geq VTC_{abl} + M \times (ZV_{abld} - 1) \qquad (44)$$

**Proof:** If $ZV_{abld} = 0$, the value of the above non-linear term is given by $VTC_{abl} \times ZV_{abld} = 0$, under any value of $VTC_{abl}$. Then under constraint (44), it follows that $ZVTC_{abld} \geq -M$ and since $ZVTC_{abld}$ has a strictly positive value, the minimizing objective function ensures that $ZVTC_{abld} = 0$. Hence, the value of the non-linear above term and $ZVTC_{abld}$ have equivalent values for this case. Also, if $ZV_{abld} = 1$, the value of the above non-linear term is given by $VTC_{abl} \times ZV_{abld} = VTC_{abl}$. Then under constraint (44), it follows that $ZVTC_{abld} \geq VTC_{abl}$ and since $ZVTC_{abld}$ has a strictly positive cost-efficiency, the minimizing objective function ensures that $ZVTC_{abld} = VTC_{abl}$.

### 4.5. Robust Fuzzy Programming Approach

An efficient and effective robust supply chain is a competitive advantage for countries and corporations, helping them cope with increasing turbulence and competitive pressure.

To develop the initial model in an uncertain environment, the demand is assumed to be uncertain to cope with uncertainty in real-life optimization problems. In this regard, the Robust Fuzzy Programming Approach introduced by Pishvaee and Fazli Khalaf (35) is adopted to handle uncertainty.

Constraint (41) is the only constraint that contains uncertain parameters. So, according to (35), this constraint is modified as follows:

$$Demand_{dpa} - (t^m + \frac{\varphi_t - \varphi_t'}{3})(1-\alpha) \leq \sum_l \sum_b UD_{bapld} * zv_{bald} + UnmD_{adp} \qquad \forall a \neq 1, d, p, b \neq a, 0 \leq \alpha \leq 1 \qquad (45)$$

According to Cadenas and Verdegay (36) and Piedro et al. (37), one fuzzy number $\tilde{t}$ is used to show the v action of the soft constraint.

Now with the assumption of $\tilde{t}$ being a triangular fuzzy number, based on the Yager (38), it can be written as $(t^m + \frac{\varphi_t - \varphi_t'}{3})$, and the parameter $\alpha$ is an indicator of the minimum satisfaction level of the flexible constraint. So, the term is the possible violation of the constraint. With this change, the new cost term $\gamma[(t^m + \frac{\varphi_t - \varphi_t'}{3})(1-\alpha)]$ will be added to the objective function, meaning the total cost of possible violation on the soft constraint where $\gamma$ is the penalty value. For more information, see the [38].

## 5. Solution Method

This section presents a heuristic method based on the genetic algorithm to solve the problem.

### 5.1. Genetic Algorithm

GA algorithm was initially studied by John Holland (39). GA is used to optimize different scheduling problems. Genetic algorithms are general-purpose stochastic search procedures that mimic the process of natural selection. A population is a set of possible solutions or chromosomes that have undergone reproduction (selection and modification) to form the next generation of solutions (offspring). The objective function evaluates the quality of solutions. This process continues until the stopping criterion is reached. Genetic algorithms are discussed in more detail in Holland and Goldberg. Figure 2 shows the procedures of the proposed GA algorithm.

1. **Input**: $npop, pcross, pmut, Miter$
2. **Output**: "$\delta_{Best}, obj(\delta_{Best})$"
3. **Begin**
4. Randomly generate the initial solutions with the heuristic method and then use heuristic1 and heuristic2 to calculate the objective function for each initial solution and select the best objective.

7. $iter = 1$
8. **While** $iter < M_{iter}$

```
            for i = 1 to cross (pcross*npop)
                Select two solutions with the roulette-wheel selection method
                if iter < 0.33*M_iter
                    Generate two new solutions based on heuristic methods with the use of chromosomes 1,2
                end
                if iter > 0.33*M_iter and iter < 0.66*M_iter
                    Generate two new solutions based on heuristic methods on chromosome 1
                end
                if iter > 0.66*M_iter
                    Generate two new solutions based on heuristic methods on chromosome 2
                end
            end
            for i = 1 to nmut (pmut*npop)
                Select two solutions with the roulette-wheel selection method
                if iter < 0.8*M_iter
                    Generate a new solution based on swap operation with the use of chromosomes 1,2
                end
                if iter > 0.8*M_iter
                    Generate a new solution based on swap operation with the use of chromosome 1
                end
            end
            calculate objective function for new solution based on mut and cross
            select the best solutions (the number of pop solutions selected)

23.         iter = iter+1
24.     return "
    25. end
```

**Figure 2.** The procedure of the proposed GA-based heuristic

### 5.1.1. Chromosomes

In the genetic algorithm, each chromosome represents a point in the search space and a possible solution to the problem. Two strings are used to show the chromosome. The first string shows the production sequence of product families in the production lines. It starts production in terms of this sequence. Whenever a period ends, it stops production and starts the production of the following family in the next period. The number of rows at the string equals the number of production lines. The string length is three times the total number of production periods (3D). Such a considerable length can remove any production constraint, such as the allowable limit on producing a product family in a line. Figure 3 shows an example of two production families, three time periods, and two production lines.

| Line 1 | 1 | 1 | 2 | 1 | 2 | 2 | 1 | 2 | 2 |
| Line 2 | 1 | 1 | 1 | 2 | 2 | 1 | 1 | 1 | 2 |

**Figure 3.** String of line

The second string shows the delivery sequence of products. It delivers in this sequence, considering distribution centers and the number and capacity of vehicles. The number of rows at the string equals the number of cars. Figure 4 shows an example of three vehicles and ten distribution centers.

| Vehicle 1 | 5  | 9 | 10 | 8 | 2 | 7 | 4 | 3 | 1 | 6 |
| Vehicle 2 | 10 | 1 | 7  | 5 | 4 | 3 | 9 | 2 | 6 | 8 |
| Vehicle 3 | 5  | 9 | 10 | 1 | 8 | 4 | 3 | 6 | 7 | 2 |

**Figure 4.** String of vehicle

### 5.1.2. Initial Population

The number of the initial population in the genetic algorithm is equal to 100. To generate the initial population, we perform the following procedure
- The first string

As a first step, each product family's preparation and processing costs are set in ascending order for each line. Then, every one of the first three families in each line is randomly distributed across all the days needed for that line. If no product family is produced, that family is randomly distributed across the production lines. It provides an initial solution to the problem. Forty-nine additional solutions are obtained by randomly displacing product families in different lines. In this way, two lines are randomly selected from product families and are replaced with each other at a distance of the chromosome length in two lines.
- The second chromosome

The second string representing the delivery sequence of products to distribution centers for each vehicle is randomly generated. In this case, one sequence is allocated to each vehicle for the delivery of products.

### 5.1.3. Fitness Value

In this section, the objective function of the problem is calculated. The first string determines how to produce product families in the production lines. The product sequence is specified by calculating the changeover cost for all products belonging to the same family in a production line. Then, the following two methods are applied to determine the number of productions per line.

**First Method:** Production based on demand for each day and production line (*PBD*)
Considering the cooling process for which a certain period is required, the number of products to be produced is precise. In this situation, the number of productions on each day is equal to the next day's demand $P$, and $P$ is considered equal to the time required for cooling the products. Then lines of the problem would be arranged randomly. Considering the first line, production will begin according to the first string on that line. Here, the number of products would appear under the standard period and equal to the next day's demand $P$ Then consider the second line; its products will also be produced. If the number of products required to be produced on that day reaches its quorum, no more will be produced from that product. This will continue until all production required for one day is produced and/or based on the capacity of production lines. After a review of all production lines, if the number of products required to be produced has not reached the number demanded, i.e., there would be a shortage, this action will be continued until all lines reach the number of required products for that day If production capacity is more than the demand of the day; the empty capacity would be stored in a separate matrix; and if a shortage occurs in the following periods and we would be able to produce that product on that line and the same day; then, it will be produced (such restrictions as deterioration of products and shelf life would also be studied under such situation so that no product would be available in the system after the end of its shelf life). This process will be done for all days so that all required products are produced.
Clarifying production at the end of each day would be sent respectively to be cooled and dispatched. Here, products would be distributed based on the second chromosome. As far as the sequence path is clear for each vehicle, the cost of dispatch of products would be calculated for all paths in the sequence. The minimum cost would be considered, and the same vehicle would be used to start product distribution. This will continue until all products are delivered to all distribution centers.

**Second Method:** Production based on the average production amount required for each day (*PBA*)
Similar to the previous method, the amount of product required to be produced each day will be calculated. Then, considering the number of production lines and dividing the demand for a product in one period by the number of the lines, the production amount required to be produced in one line and during one period would be clarified, and production would begin. If the number of solutions required on each line reached its required value, production would be stopped; otherwise, it would be continued until regular time and overtime production hours would be ended. Here, producing this category will end, and the following category will be produced based on the first chromosome during the next period. After the end of the production and cooling process, the product's delivery stage would be considered. Considering the demand for each product and parameter in each distribution center, those with the highest shortage cost will be considered and allocated with the highest demand. So, a few would be selected if the number of parameters is equal among different centers based on transportation costs. Products would be delivered to distribution centers by considering the vehicle's capacity and product delivery sequence in the second chromosome. The final objective function equals the best amount of the *PBA* and *PBD* methods.

### 5.1.4. Selection

The proposed GA selects a new population through the well-known roulette-wheel selection method for the next generation. This method selects a new population based on the probability distribution associated with the fitness values of chromosomes.

### 5.1.5. Crossover

After the selection phase and choosing the chromosomes, a common point of these strings would be selected randomly for the first string of each chromosome (sequence of product categories in production lines). Then, the crossover would be applied to it. An example of this process is shown in Figure 5.

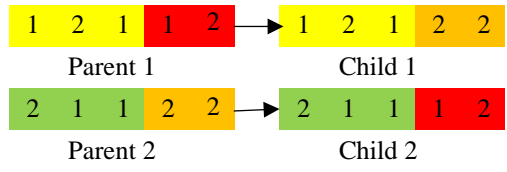
**Figure 5.** Example of crossover operation

A single-point crossover has also been used for crossover on the second string of each chromosome. The critical point in terms of the second chromosome is the lack of repetitive values placed on each chromosome strand which has to include all distribution centers so that it would not be received by any of the distribution centers twice. For string modification and considering production cost with the previous state of it, the order of placement of distribution centers in chromosome would be considered so that transportation cost would not be increased compared to the previous state for the same chromosome. Since transportation cost is precise, the repetitive state in chromosomes would be omitted so that no cost would be added. An example of this process is shown in Fig. 6.

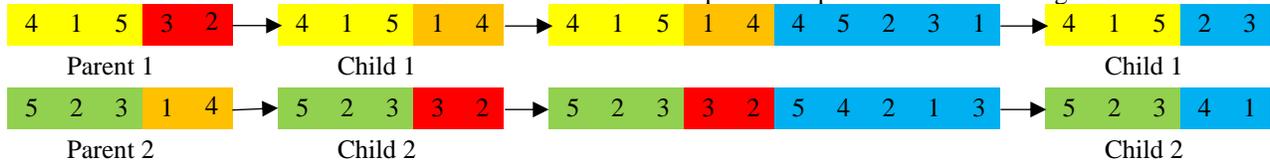
**Figure 6.** Example of the crossover operation

By considering the number of repetitions in the problem, different crossovers would be applied.

### 5.1.6. Mutation
For mutation, one chromosome will be selected randomly. Here, the swap has been used for mutation. In the selected chromosome, primarily, two lines would be selected randomly for the first part of the string, which is the same production sequence of the product category in production lines. Then, one gene from each selected line will be chosen randomly, and their locations will be replaced. If the two selected genes from the two concerned lines are of one product category, a swap will not occur, and two other genes will be selected. For mutation on the second chromosome, two genes will be randomly selected and swapped with each other. Figure7 shows an example of the mutation operation.

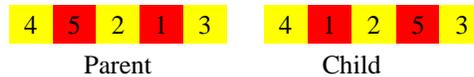
**Figure7.** Example of mutation operation

### 5.1.8. Reproduction
Crossover and mutation (genetic operators) probabilities determine the rate of reproduction. When the components of the GA are appropriately chosen, the reproduction process continually improves the population, finally converging on solutions close to a globally optimal solution. In each generation of the algorithm, the current population and the new population generated by genetic operators are combined. Then, the solutions' fitness value is evaluated to sort them in descending order. The population size is 100. Therefore, the 100 individuals with the best fitness value are preserved for the next generation.

### 5.1.9. Stopping Criterion
The stopping criterion combines the maximum number of generations and convergence. The algorithm stops when the maximum number of generations is 300. Another stopping criterion is no improvement in the best fitness value available after 100 generations. The algorithm is terminated when each of these stopping criteria is satisfied. From repletion 1 to 100, the crossover would be applied on the first string and mutation on the second string only. From repletion 100 to 20, quite the contrary, the crossover would be performed on the second string and mutation on the first string. Finally, from repletion 200 to 300, both the crossover and mutation would be employed on either string.

## 6. RESULTS
In this section, we first describe how to generate sample problems, and then describe the results obtained from solving the problem using the heuristic algorithm. Then a real-world case is considered, and its results are analyzed using the proposed model. First, the proposed model was linearized, and then the calculations were done using *MATLAB* and CPLEX software on a Corei2 computer at 2.4 GHZ with 2 GB of DDRII 533 MHZ RAM.

### 6.1. Data Generation

The proposed problem is investigated with two types of data. The first type is a sample problem in three different sizes, and the second type of data has been extracted from a real case study in Tehran.

### 6.1.1. Numerical Example
Case problems are in three scales, small, medium, and large. Ten problem cases were run by CPLEX and the proposed algorithm in each problem size, and the results were analyzed. Table 3 shows the sets and parameter generation scheme.

**TABLE 3.** Test problem and Numerical example scheme.

| Sets | Small | Medium | Large | Sets | Small | Medium | Large |
|---|---|---|---|---|---|---|---|
| $a$ | Uniform (3,5) | Uniform (6,10) | Uniform (11,15) | $p$ | Uniform (2,4) | Uniform (5,8) | Uniform (8,12) |
| $f$ | Uniform (1,2) | Uniform (3,4) | Uniform (5,6) | $d$ | Uniform (1,3) | Uniform (4,7) | Uniform (8,12) |
| $j$ | Uniform (2,3) | Uniform (4,5) | Uniform (6,8) | $r$ | Uniform (1,2) | Uniform (3,4) | Uniform (5,6) |
| $l$ | Uniform (2,3) | Uniform (4,5) | Uniform (6,7) | $i$ | Uniform (1,3) | Uniform (4,7) | Uniform (8,12) |

| parameter | Generation scheme | parameter | Generation scheme | parameter | Generation scheme | parameter | Generation scheme |
|---|---|---|---|---|---|---|---|
| $\Pr et_r$ | Uniform (3,5) | $MaxTC_l$ | Uniform (1000,2000) | $Setup_{jp}$ | Uniform (30,50) | $MaxLots_p$ | Uniform (3000,5000) |
| $Rtime_i$ | Uniform (420,480) | $MinTC_l$ | Uniform (300,600) | $Relt_{ri}$ | Uniform (5,15) | $MinLots_p$ | Uniform (500,2500) |
| $Maxtime_i$ | Uniform (540,720) | $Pcapacity_i$ | Uniform (500,5000) | $O_{ji}$ | Uniform (5,10) | $Demand_{dpa}$ | Uniform (50,100) |
| $ShelfLife_p$ | Uniform (15,45) | $Pallet_p$ | Uniform (0.05,0.08) | $Bpc_{ri}$ | Uniform (5,10) | $Prate_{jp}$ | Uniform (5,10) |
| $CrRate_p$ | Random (0,1) | $dailyop_{ji}$ | Uniform (5,10) | $IC_{pi}$ | Uniform (0.05,2) | $Unmd\,Cost_{ap}$ | Uniform (5,15) |
| $FCT_l$ | Uniform (100,150) | $dailysh_{ji}$ | Uniform (10,20) | $VTC_{abl}$ | Uniform (1,10) | $F\,Cost_{ji}$ | Uniform (15,30) |
| $Var\,Cost_p$ | Uniform (0.09,0.15) | $W_{ji}$ | Uniform (480,720) | $Chc_{feji}$ | Uniform (5,10) | $MuMax_{ri}$ | Uniform (2000,4000) |
| $Overt\,Cost_i$ | Uniform (1,3) | $Cht_{fej}$ | Uniform (5,10) | $Line\,Cost_{ji}$ | Uniform (10,15) | $MuMin_{ri}$ | Uniform (300,2000) |
| $CST$ | Uniform (420,480) | $QCTime$ | Uniform (240,360) | $StCapacity$ | Uniform (500,5000) | | |

### 6.2. Numerical Results
In this section, the results of the proposed method are examined with simulated data. The following table compares the proposed algorithm with the optimal solution obtained by the branch and cut algorithm via the CPLEX solver. In this case, the deviation of the proposed method from the exact answer is calculated as the $\frac{Z_{GA}-Z_{Cplex}}{Z_{GA}}$.

Note that we interrupt the solver after 1h (3600 seconds) and the objective bounds are reported as the lower bound if the software doesn't reach the optimal solution (for medium and large-size problems).

**TABLE 4.** Comparison between CPLEX and Heuristic method in small size problem

| Sample size | DC | f | j | l | p | d | r | i | %Gap between CPLEX & Heuristic method | Time (s) CPLEX | Time (s) Heuristic method |
|---|---|---|---|---|---|---|---|---|---|---|---|
| 1 | 3 | 1 | 1 | 2 | 2 | 1 | 1 | 1 | 0 | 36 | 12.6 |
| 2 | 3 | 1 | 2 | 2 | 2 | 1 | 1 | 1 | 1.3 | 570 | 21.1 |
| 3 | 3 | 1 | 2 | 2 | 2 | 2 | 1 | 2 | 0 | 2163 | 18.6 |
| 4 | 4 | 2 | 1 | 2 | 3 | 2 | 2 | 2 | 0 | 117 | 17.3 |
| 5 | 4 | 2 | 1 | 2 | 3 | 2 | 2 | 2 | 0 | 1861 | 21.4 |
| 6 | 4 | 2 | 2 | 2 | 3 | 3 | 2 | 3 | 1.1 | 2968 | 22.3 |
| 7 | 5 | 2 | 1 | 2 | 4 | 3 | 2 | 3 | 0.6 | 1630 | 16.1 |
| 8 | 5 | 2 | 2 | 3 | 4 | 3 | 2 | 3 | 1.3 | 1714 | 21.5 |
| 9 | 5 | 2 | 1 | 3 | 4 | 3 | 2 | 3 | 0 | 2980 | 11.3 |
| 10 | 5 | 2 | 2 | 3 | 4 | 3 | 2 | 3 | 0 | 2630 | 10.9 |

It can be seen from Table 4 that the proposed method was able to reach an optimum solution of 6 examples out of 9. The mean deviation of the proposed method from the exact answer in all examples in the table above is 0.4% which indicates the efficiency of the proposed method for solving problems. The average solution time of the CPLEX method is 1667 seconds, and the proposed algorithm is 17.31 seconds.

### 6.3. Sensitivity Analysis

In the following tables, we investigate the changes in *QB*, *II*, and Unmet Demand quantities according to different amounts of *CrRate*, Shelf life, and Satisfaction level. In these charts, *QB* means the average output per day, *II* is the average quantity of products available per day, and Unmet Demand is the average quantity of products in the distribution centers. Also, the average objective function improvement in the proposed algorithm compared to the CPLEX method is computed by the $\frac{Z_{Cplex} - Z_{GA}}{Z_{Cplex}}$.

**TABLE 5.** compares between CPLEX and proposed GA in medium size under different CrRate

| $CrRate_p$ | CPLEX | | | | | GA | | | | |
|---|---|---|---|---|---|---|---|---|---|---|
| | z | $QB_{pi}$ | $ll_{pi}$ | $UnmD_{adp}$ | Time | z | $QB_{pi}$ | $ll_{pi}$ | $UnmD_{adp}$ | Time |
| 0.1 | 214,010 | 375 | 139 | 80 | 3600 | 192,100 | 409 | 131 | 73 | 46.59 |
| 0.2 | 220,600 | 364 | 144 | 83 | 3600 | 194,128 | 400 | 134 | 78 | 47.59 |
| 0.3 | 283,981 | 348 | 139 | 90 | 3600 | 258,423 | 382 | 125 | 84 | 46.98 |
| 0.4 | 345,123 | 341 | 132 | 99 | 3600 | 307,159 | 371 | 120 | 90 | 48.61 |
| 0.5 | 393,458 | 332 | 127 | 106 | 3600 | 346,243 | 357 | 114 | 100 | 51.65 |
| 0.6 | 458,968 | 311 | 121 | 113 | 3600 | 426,840 | 346 | 108 | 100 | 52.27 |
| 0.7 | 510,315 | 290 | 108 | 121 | 3600 | 449,077 | 330 | 102 | 110 | 54.54 |
| 0.8 | 564,568 | 264 | 100 | 129 | 3600 | 519,403 | 300 | 89 | 119 | 57.65 |
| 0.9 | 595,710 | 251 | 104 | 134 | 3600 | 524,300 | 286 | 94 | 124 | 58.61 |

As seen from Table 5, the *QB* and *II* when *CrRate* is equal to 0.1 compared to the time when *CrRate* is 0.9 in the CPLEX method, and the proposed algorithm has seen a 33%, 30%, 25%, and 28% fall, respectively. On the other hand, the average amount of Unmet Demand increased by 40% and 41%, respectively. Increasing the amount of *CrRate* means that the customer tends to buy fresher products with a longer shelf life. So, in this case, the products are produced as close to the demand day as possible, and due to capacity constraints, the amount of Unmet Demand increases under heavy demand. The *QB*, *II*, and Unmet Demand in the proposed algorithm are 12%, 8.7%, and 8% higher than the CPLEX method. The average objective function has seen a 10% improvement for all *CrRates* in the proposed model compared to CPLEX as well. The solution time in the CPLEX method is 3600 seconds, and the average solution time of the proposed algorithm is 51.61.

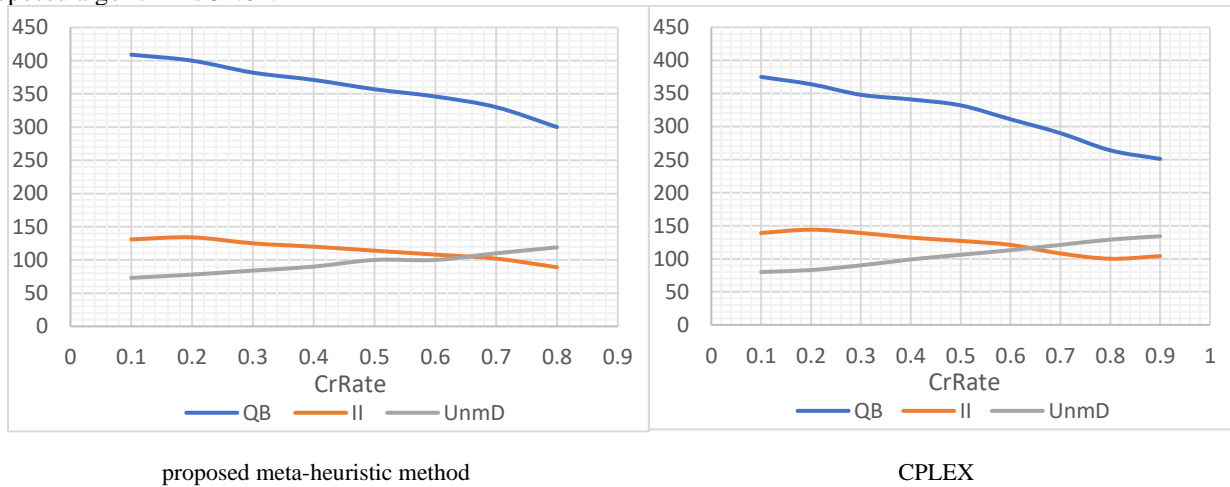

proposed meta-heuristic method            CPLEX

**Figure 8.** Changes in production amount, inventory level, and Unmet demand under different CrRate levels.

According to Table 5 and Figure 8, the production policy of the company transforms in a way that, daily, they produce precisely the amount of demand. The company is shifting to producing smaller batches rather than producing larger ones. Moreover, the smaller the batch size, the more sterilization and cleaning operations are performed, resulting in higher production costs. On the other hand, the inventory level is reduced under these conditions. Thus, although the maintenance costs decrease, total costs increase.

**TABLE 6.** Comparing CPLEX and proposed GA in large size under different Shelf life

| $Shelflife_p$ | CPLEX | | | | | GA | | | | |
|---|---|---|---|---|---|---|---|---|---|---|
| | Z | $QB_{pi}$ | $II_{pi}$ | $UnmD_{adp}$ | Time(s) | Z | $QB_{pi}$ | $II_{pi}$ | $UnmD_{adp}$ | Time(s) |
| 15 | 220,600 | 198 | 48 | 46 | 3600 | 198,540 | 215 | 45 | 43 | 48.35 |
| 20 | 201,236 | 315 | 62 | 43 | 3600 | 189,162 | 346 | 56 | 41 | 48.36 |
| 25 | 183,254 | 386 | 75 | 40 | 3600 | 166,761 | 429 | 68 | 36 | 49.71 |
| 30 | 152,345 | 481 | 87 | 37 | 3600 | 140,157 | 512 | 80 | 35 | 52.48 |
| 35 | 124,562 | 551 | 97 | 34 | 3600 | 114,597 | 599 | 92 | 32 | 53.68 |
| 40 | 94,263 | 624 | 109 | 30 | 3600 | 89,550 | 678 | 98 | 28 | 55.19 |
| 45 | 53,245 | 697 | 122 | 23 | 3600 | 50,050 | 741 | 114 | 21 | 58.87 |

In Table 6, by increasing the Shelf life, the products could be produced earlier and maintained longer, which reduces Unmet Demand. $QB$ and $II$ per day when Shelf life is equal to 15 compared to the time when it is 45 in the CPLEX method, and the proposed algorithm has seen 240%, 250% and, 154%, 153% rise, respectively. On the other hand, the amount of Unmet Demand decreased by 66% and 67%, respectively. As a result, by increasing the Shelf life, $QB$ and $II$ increase, and Unmet Demand decreases. Unmet Demand and objective function in the proposed algorithm are 7%, 8.5%, 7.8%, and 6.7% less than those in CPLEX, respectively. The average solution time of the proposed algorithm is 52.37.

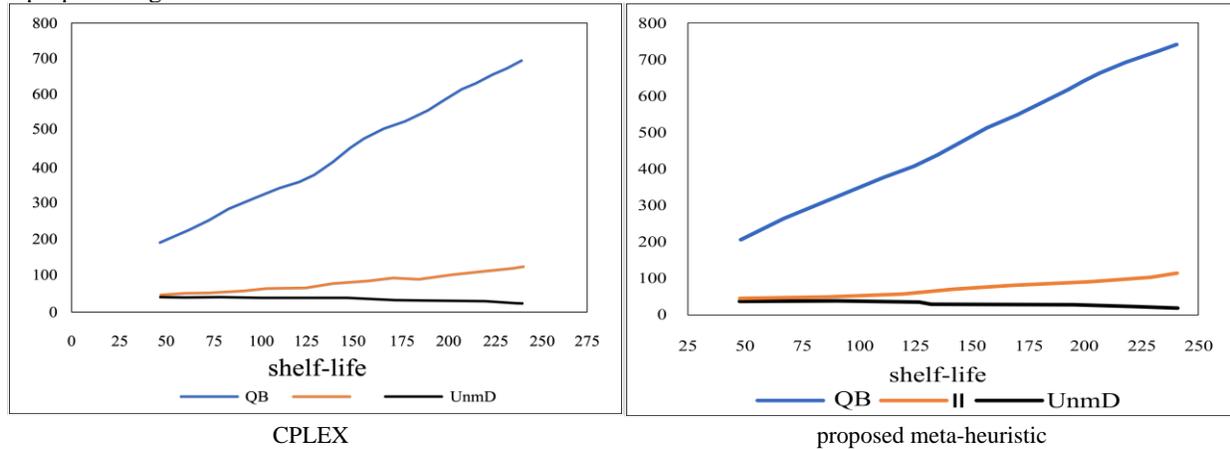

CPLEX                proposed meta-heuristic

**Figure 9.** Changes in production amount, inventory level, and Unmet demand under different Shelf-life levels

With a longer Shelf life, the manufacturer can allocate more of its stock to storing products so that it does not face a shortage in a highly uncertain situation. In this case, sterilization and cleaning operations will be reduced as the production cost (see Figure 9).

**TABLE 7.** Comparing CPLEX and proposed GA in large size under different levels of satisfaction

| $\alpha$ | CPLEX | | | | | GA | | | | |
|---|---|---|---|---|---|---|---|---|---|---|
| | Z | $QB_{pi}$ | $II_{pi}$ | $UnmD_{adp}$ | Time(s) | Z | $QB_{pi}$ | $II_{pi}$ | $UnmD_{adp}$ | Time(s) |
| 1 | 198,650 | 272 | 45 | 71 | 3600 | 182,758 | 296 | 41 | 65 | 50.35 |
| 0.8 | 220,600 | 279 | 48 | 74 | 3600 | 200,746 | 310 | 44 | 68 | 52.36 |
| 0.7 | 275,469 | 300 | 62 | 66 | 3600 | 253,431 | 330 | 57 | 60 | 51.71 |
| 0.6 | 336,241 | 311 | 75 | 58 | 3600 | 302,617 | 346 | 68 | 53 | 54.48 |
| 0.5 | 384,526 | 321 | 85 | 50 | 3600 | 349,919 | 357 | 77 | 46 | 56.68 |
| 0.4 | 457,846 | 334 | 97 | 41 | 3600 | 412,061 | 371 | 88 | 37 | 58.19 |
| 0.3 | 508,562 | 351 | 109 | 32 | 3600 | 462,791 | 382 | 101 | 29 | 60.87 |
| 0.2 | 565,821 | 351 | 109 | 43 | 3600 | 509,239 | 382 | 101 | 41 | 62.95 |
| 0.1 | 597,364 | 351 | 109 | 57 | 3600 | 543,601 | 382 | 101 | 59 | 63.78 |

As can be seen from Table 7, the level of $QB$ and $II$ when the satisfaction level is 1, compared to the time when it is equal to 0.1 in the CPLEX method and the proposed algorithm have increased 29.04%, 29.05%, and 142%, 146% respectively. In contrast, Unmet Demand decreased by 19.7% and 9.2%, respectively. *QB*, *II* Unmet Demand, and objective function in the proposed algorithm are 6.9%, 8.3%, 9.9%, and 9.1% higher than the CPLEX method. The average solution time of the algorithm presented in the diagrams below shows the number of changes in *QB*, *II,* and Unmet Demand versus the level of uncertainty.

Decreasing the level of satisfaction increases the level of *QB* and *II,* reducing the amount of Unmet Demand. However, in the event of severe fluctuations in the business environment, the company's production policies will move in a way that can cope with the fluctuations, which of course, require high costs such as inventory maintenance costs. However, since the level of satisfaction reaches 0.3, the amount of *QB* and *II* remains constant due to production and inventory constraints. So, in this case, the amount of Unmet Demand rises.

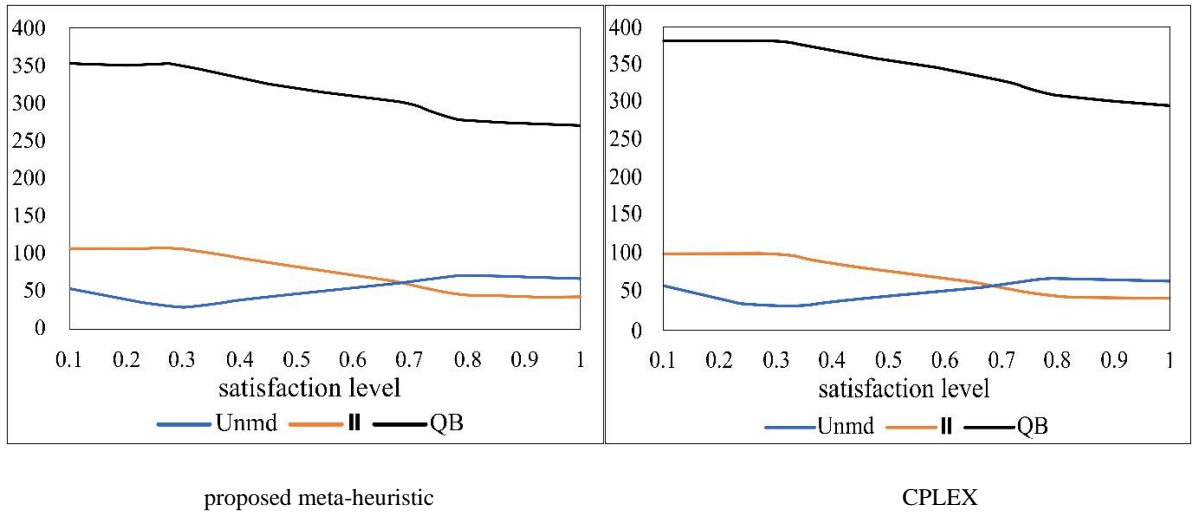

proposed meta-heuristic            CPLEX
**Figure 10.** Changes in production amount, inventory level, and Unmet demand under different Uncertainty levels

In summary, with increasing levels of uncertainty and fluctuations in demand, inventory levels increase which leads to a massive rise in inventory costs and subsequently exacerbates the challenge of keeping products fresh (see Figure 10).

### 6.4. Numerical solution

The real case problem is based on the supply chain of a dairy company in Iran, and the structure of the supply chain network is made according to the study of this company. In the case study, a single plant delivers the final products to 10 Canbo branches (DCs). The plant includes six different yogurt products. Each product has specific features. The products with similar characteristics are placed into product families (see Table 8).

**TABLE 8**. Allocation of each product to the family of products

| Family | Set Yogurt | Stirred Yogurt |
|---|---|---|
| Product | 1. Cream Yogurt | 1. Eggplant Yogurt |
|  | 2. Low-fat Yogurt | 2. Strawberry Fruit Yogurt |
|  | 3. Traditional strained Yogurt | 3. Cucumber Yogurt |

The planning and scheduling are performed in a short-term horizon of one week (five working days). There are two identical production lines, each of which can produce all products. There are 8 hours available for regular working time daily. Under heavy demand, overtime of about 4 hours is allowed every working day. After finishing the production process, to achieve ultimate stability and preserve the high quality of the final product, they place packed products in cooling storage containers for one day. Quality control is also addressed in this phase. The products are allowed to be delivered after they complete this stage. The production site delivers products to distribution centers with 30 items in each package. The plant has contracts with three vehicles that deliver the final products to the distribution centers. Each vehicle is characterized by a Min and Max capacity of 500 and 1500 Yogurt products, respectively. The proposed model has been applied to this case study. The amount of demand for each Canbo branch for all kinds of products in 5 periods is demonstrated in Table 9 (Appendix A). The amount of production on each production line and each period is also shown in Table 10.

**TABLE 10.** Computational result for the case study

| Production Line | Yogurt Products | Period 1 | | Period 2 | | Period 3 | | Period 4 | | Period 5 | |
|---|---|---|---|---|---|---|---|---|---|---|---|
| | | Model | Real | Model | Real | Model | Real | Model | Real | Model | Real |
| Line 1 | Cream yogurt | 430 | 430 | 380 | 380 | 375 | 375 | 350 | 350 | 355 | 355 |
| | Strawberry fruit yogurt | 321 | 321 | 331 | 331 | 324 | 324 | 329 | 329 | 317 | 317 |
| | Traditional strained yogurt | 360 | 360 | 345 | 345 | 355 | 355 | 355 | 355 | 350 | 350 |
| Line 2 | Eggplant yogurt | 327 | 327 | 322 | 322 | 367 | 367 | 342 | 342 | 337 | 337 |
| | Low-fat yogurt | 395 | 395 | 365 | 365 | 365 | 365 | 370 | 370 | 395 | 395 |
| | Cucumber yogurt | 312 | 312 | 307 | 307 | 312 | 312 | 321 | 321 | 318 | 318 |

**TABLE 11.** Detail results of the routing problem

| period | | Vehicle schedule | Number of served DCs | Total amount of transported products |
|---|---|---|---|---|
| Period 1 | Route 1 | 1 → 2 → 3 → 4 | 4 | 1417 |
| | Route 2 | 5 → 6 → 7 → 8 | 4 | 1497 |
| | Route 3 | 9 → 10 | 2 | 1015 |
| Period 2 | Route 1 | 1 → 2 → 3 | 3 | 1172 |
| | Route 2 | 4 → 5 → 6 → 7 | 4 | 1346 |
| | Route 3 | 8 → 9 → 10 | 3 | 1007 |
| Period 3 | Route 1 | 1 → 2 → 3 | 3 | 1187 |
| | Route 2 | 4 → 5 → 6 → 7 | 4 | 1368 |
| | Route 3 | 8 → 9 → 10 | 3 | 997 |
| Period 4 | Route 1 | 1 → 2 → 3 | 3 | 1214 |
| | Route 2 | 4 → 5 → 6 → 7 → 8 | 5 | 1498 |
| | Route 3 | 9 → 10 | 2 | 900 |
| Period 5 | Route 1 | 1 → 2 → 3 | 3 | 1297 |
| | Route 2 | 4 → 5 → 6 | 3 | 1287 |
| | Route 3 | 7 → 8 → 9 → 10 | 4 | 1484 |

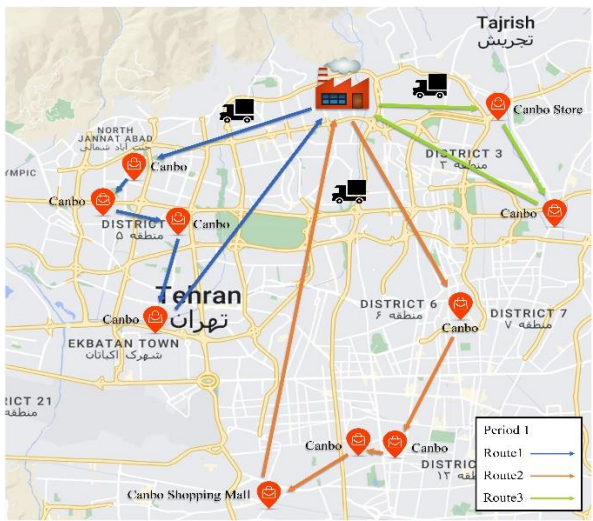
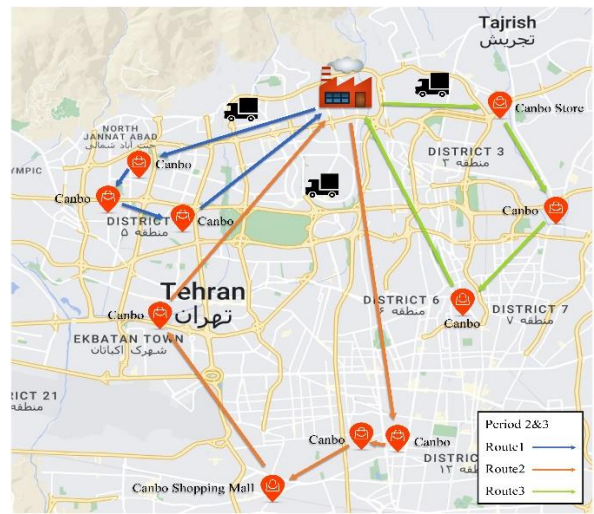

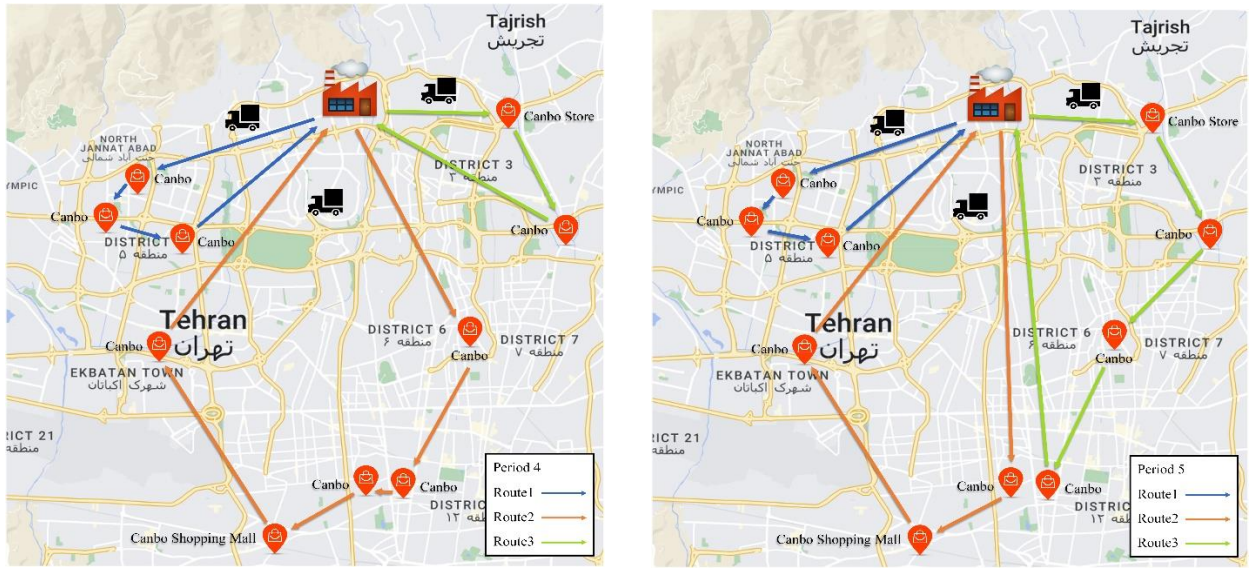

**Figure11.** Geographical description of the operation route of each vehicle in each period

The above figures show how products are delivered to each DC on the Tehran city map. In this case, the way that vehicles are assigned to each DC varies depending on the different demands in different periods. Table 11 illustrates detailed delivery results of vehicles, including specific routes and the total amount of products transported by each vehicle. Because each vehicle is characterized by a Max and Min capacity of 500 and 1500, respectively, they cannot carry more than 1500 products in one run due to capacity constraints. For example, as shown in Figure 11, in period 1, vehicle 1 delivers the demand of DCs 1, 2, 3, and 4 in the first route. But according to Table 9, there has been an increase in the demand in the following periods, so vehicle 1 cannot carry the demand for DCs 1-4 anymore (since the total demand is more than Max capacity of the vehicle, see Table 11 and Table 9), so the Canbo number 4 is served by vehicle 2.

## 7. CONCLUSION

In this study, we investigate the problem of production and distribution planning, scheduling, and routing inspired by studying a dairy industry supply chain network. The considered problem is a multi-product, single plant, multi-distribution center, multi-period, and multi-transportation track. This is followed by a decision on the allocation of vehicles, which are characterized by Max and Min capacity, to the distribution centers, inventory and production level of each product in each period, and the required amount of overtime under the condition of heavy demand. This research proposes a mixed-integer linear programming (MILP) model to minimize the total costs, including production, setup, overtime, unmet demand, and transportation. The robust fuzzy optimization approach is employed to deal with the uncertainty of the demand product in each distribution center. Moreover, to demonstrate the applicability of the proposed model, we provide a real-world case study from the Kalleh dairy industry in Iran. The proposed genetic algorithm (GA) is applied to solve the extended problem. The obtained results show a relative gap of about 0.4% between GA and the optimum solution in terms of the objective function value found by the CPLEX solver (Branch and cut algorithm) in small-size problems. The performance of GA is evaluated and compared by medium and large-sized problems concerning the some of defined measures. The experimental results confirm that increasing the shelf-life parameter leads to a rise in stocking products to meet unanticipated demand. Therefore, production costs are reduced, and the total costs are reduced naturally. The proposed model considers uncertainty parameters and indicates that with increasing levels of uncertainty and fluctuations in demand, inventory levels increase, leading to a fall in unmet demand but a massive rise in inventory costs and exacerbating the challenge of keeping products fresh. Considering the sustainability dimensions in the proposed mathematical model is recommended for future research. Furthermore, proposing another solution approach is an interesting research field to enhance solution performance.

**Funding:** This research received no external funding
**Data Availability Statement:** The data presented in this study are available on request from the corresponding author.
**Conflict of interest:** The authors declare no conflicts of interest.

# Appendix A

**TABLE 9.** The demand of each distribution center per product

| period | Product Type | Distribution centers | | | | | | | | | |
|---|---|---|---|---|---|---|---|---|---|---|---|
| | | Canbo 1 | Canbo 2 | Canbo 3 | Canbo 4 | Canbo 5 | Canbo 6 | Canbo 7 | Canbo 8 | Canbo 9 | Canbo 10 |
| Period 1 | Cream yogurt | 65 | 60 | 65 | 65 | 60 | 60 | 55 | 95 | 90 | 85 |
| | Eggplant yogurt | 50 | 50 | 50 | 50 | 50 | 42 | 55 | 87 | 85 | 80 |
| | Low-fat yogurt | 60 | 65 | 60 | 70 | 65 | 55 | 60 | 85 | 90 | 85 |
| | Strawberry fruit yogurt | 65 | 52 | 50 | 55 | 50 | 50 | 52 | 82 | 85 | 85 |
| | Traditional strained yogurt | 60 | 65 | 75 | 65 | 55 | 60 | 50 | 92 | 85 | 85 |
| | Cucumber yogurt | 50 | 50 | 65 | 55 | 52 | 50 | 55 | 80 | 80 | 80 |
| | Total demand | 350 | 342 | 365 | 360 | 332 | 317 | 327 | 521 | 515 | 500 |
| Period 2 | Cream yogurt | 65 | 75 | 65 | 75 | 60 | 60 | 55 | 65 | 60 | 55 |
| | Eggplant yogurt | 60 | 65 | 60 | 80 | 55 | 52 | 55 | 55 | 60 | 50 |
| | Low-fat yogurt | 85 | 65 | 70 | 65 | 65 | 55 | 60 | 55 | 70 | 55 |
| | Strawberry fruit yogurt | 65 | 72 | 65 | 55 | 50 | 50 | 52 | 52 | 55 | 55 |
| | Traditional strained yogurt | 75 | 60 | 75 | 60 | 55 | 60 | 50 | 60 | 55 | 55 |
| | Cucumber yogurt | 55 | 65 | 50 | 50 | 52 | 50 | 55 | 50 | 50 | 50 |
| | Total demand | 405 | 402 | 385 | 355 | 337 | 327 | 327 | 337 | 350 | 320 |
| Period 3 | Cream yogurt | 85 | 65 | 75 | 60 | 60 | 60 | 55 | 60 | 65 | 55 |
| | Eggplant yogurt | 75 | 50 | 50 | 60 | 55 | 52 | 55 | 55 | 55 | 50 |
| | Low-fat yogurt | 65 | 85 | 80 | 75 | 65 | 55 | 60 | 55 | 60 | 55 |
| | Strawberry fruit yogurt | 60 | 55 | 50 | 55 | 50 | 50 | 52 | 52 | 55 | 55 |
| | Traditional strained yogurt | 82 | 80 | 75 | 82 | 55 | 55 | 50 | 60 | 55 | 55 |
| | Cucumber yogurt | 50 | 60 | 50 | 50 | 52 | 50 | 55 | 50 | 50 | 55 |
| | Total demand | 412 | 395 | 380 | 382 | 337 | 322 | 327 | 332 | 340 | 325 |
| Period 4 | Cream yogurt | 80 | 70 | 65 | 55 | 50 | 60 | 55 | 55 | 70 | 75 |
| | Eggplant yogurt | 65 | 60 | 65 | 45 | 45 | 42 | 40 | 50 | 55 | 65 |
| | Low-fat yogurt | 70 | 75 | 60 | 55 | 50 | 50 | 40 | 55 | 90 | 85 |
| | Strawberry fruit yogurt | 55 | 60 | 72 | 50 | 45 | 45 | 50 | 52 | 80 | 70 |
| | Traditional strained yogurt | 85 | 85 | 60 | 55 | 55 | 55 | 55 | 60 | 85 | 95 |
| | Cucumber yogurt | 62 | 70 | 55 | 42 | 40 | 50 | 52 | 45 | 60 | 70 |
| | Total demand | 417 | 420 | 377 | 302 | 285 | 302 | 292 | 317 | 440 | 460 |
| Period 5 | Cream yogurt | 60 | 85 | 85 | 60 | 60 | 75 | 65 | 60 | 60 | 65 |
| | Eggplant yogurt | 85 | 60 | 85 | 60 | 70 | 65 | 55 | 65 | 50 | 60 |
| | Low-fat yogurt | 75 | 75 | 60 | 70 | 75 | 60 | 65 | 60 | 60 | 55 |
| | Strawberry fruit yogurt | 85 | 60 | 60 | 85 | 85 | 85 | 70 | 55 | 65 | 62 |
| | Traditional strained yogurt | 60 | 75 | 60 | 60 | 85 | 60 | 60 | 60 | 75 | 75 |
| | Cucumber yogurt | 80 | 62 | 85 | 85 | 62 | 85 | 55 | 62 | 55 | 60 |
| | Total demand | 445 | 417 | 435 | 420 | 437 | 430 | 370 | 362 | 375 | 377 |